\documentclass[12pt]{article}

\newtheorem{Thm}{Theorem}

\newtheorem{Lem}{Lemma}

\newtheorem{Cor}{Corollary}

\newtheorem{Prop}{Proposition}
\newtheorem{Def}{Definition}

\newtheorem{Rem}{Remark}
\begin{document}

\title{\bf\ A note on Hermitian forms and Maschke`s Theorem\footnote{Research
supported by MURST 60\%}}
\author{ M. Boratynski\\
\date{}
\small Dipartimento di Matematica\\
\small via Orabona 4\\
\small 70125 Bari, Italy\\
\small boratyn@dm.uniba.it}
\maketitle
\begin{abstract}
The aim of this note is to show that the ''usual'' proof of the Maschke
Theorem in the case of real and complex field can not be generalized at
least to the case of a field of rational functions on a non--singular curve.
\end{abstract}

\section*{Introduction}
Let a finite group \(G\)\ act linearly on a finite--dimensional vector
space \(V\)
 over a field \(k\). It is well known (Maschke`s Theorem) that in case the
characteristic of \(k\) does not divide the order of \(G\) the representation V
is completely reducible.
If \(k\) is a real or a complex field \(V\) admits a \(G\)--invariant
hermitian form \((\ ,\ ): V\times{V}\to{k}\) which is positively defined.
Then for any \(G\)--invariant subspace \(W\) of V we obtain
 \(V=W\oplus{W}^\perp\)
where
\(W^\perp\) denotes the orthogonal complement of \(W\) with respect to
\((\ ,\ )\).\ The decomposition \(V=W\oplus{W}^\perp\) holds since the 
restriction of \((\ ,\ )\) to \(W\) is nondegenerate. Moreover \(W^\perp\) is
\(G\)--invariant
because  \((\ ,\ ): V\times{V}\to{k}\) is \(G\)--invariant.\par 
The natural question arises whether every finite--dimensional \(G\)-module
admits
a \(G\)--invariant hermitian form such that its restriction to any (proper)
\(G\)--invariant submodule of \(V\) is nondegenerate in case the
 characteristic of \(k\)
does not divide the order of \(G\).\ If this were true we would obtain an 
alternative proof of the Maschke Theorem.\ The aim of this note is to
show that
the
answer is negative even when the action of \(G\) is trivial.
\vskip1cm
In the sequel \(k\) will denote a field with a fixed involution
ie.\ an automorphism
 \(\lambda\to\bar{\lambda}\) of order \(\leq 2\). We shall always
  assume that the involution is not an identity or ch $k\ne 2$ otherwise.           
{V} will denote a finite--dimensional vector space over \(k\).
\begin{Def}
A map  \((\ ,\ ): V\times{V}\to{k}\) is called a hermitian form if it is
 linear with respect to the first variable and moreover
 \(\overline{(v,w)}=(w,v)\) for all \(v\),\(w\in V\)
 
\end{Def}
\begin{Rem}
Of course we immediately obtain \((v,w_{1}+w_{2})=(v,w_{1})+(v,w_{2})\)\\and
 \((v,\lambda w)=\bar{\lambda}(v,w)\) for
 \(v\),\(w\),\(w_{1}\),\(w_{2}\in V\) and \(\lambda\in k\). In case the
  involution
 is an identity automorphism then obviously the notions of hermitian and
 symmetric forms coincide.
\end{Rem}
\par The following Lemma and Proposition are of course well known;
the proofs are enclosed only for the completeness sake.

 \begin{Lem}\label{L:a}
 Let
\((\ ,\ ): V\times{V}\to{k}\) be a non--zero hermitian form. Then there exists
\(v\in V\) such that \((v,v)\ne 0\).\\
Proof: Suppose \((v,v)=0\) for all \(v\in V\).
Then
\begin{equation}\label{E:a}
0=(v+w,v+w)=(v,w)+(w,v)
\end{equation}
for all $v$, $w\in V$. If the involution of $k$ is an identity we obtain
$2(v,w)=0$. So $(v,w)=0$ for all $v$, $w\in V$ since ch $k\ne 2$
 which is impossible.\\Let $\lambda\in V$ such that
$\lambda\ne\bar\lambda$ in case the involution is not an identity.Then
\begin{equation}\label{E:b}
0=(\lambda v,w)+(w,\lambda v)=\lambda (v,w)+\bar\lambda (w,v)
\end{equation}
It follows from (\ref{E:a}) and (\ref{E:b}) that
\[
(\lambda-\bar\lambda)(v,w)=\lambda (v,w)+\bar\lambda (w,v)=0
\]
So $(v,w)=0$ for all $v$, $w\in V$ which contradicts our hypothesis.Thus in
 both cases
there exists $v\in V$ such that $(v,v)\ne 0$.
\end{Lem}
\begin{Prop}\label{P:a}
 Let
$(\ ,\ ): V\times V\to k$ be a hermitian form. Then there exist\\$e_{1},e_{2},
\dots , e_{n}\in V$ which form a basis of $V$ and such that $(e_{i},e_{j})=0$
if $i\ne j$.
\\Proof: If dim $V=1$ then every $v\in V$ such that $v\ne 0$ will do.
Let dim~$V=~n>1$. We can assume that $(\ ,\ ): V\times V\to k$ is non--zero.
By the Lemma there exists $e\in V$ such that $(e,e)\ne 0$.
Let $\langle e\rangle$ denote the subspace generated by $e$. Put
$\langle e\rangle^{\perp }=\{v\in V\mid (v,e)=0$ which is a subspace of $V$.
Let $v\in V$. We obtain
\[
(v-((v,e)/(e,e))e,e)=(v,e)-((v,e)/(e,e))(e,e)=
(v,e)-(v,e)=0
\]
So
\[
v=((v,e)/(e,e))e+[v-((v,e)/(e,e))e] 
\]
with
\[
((v,e)/(e,e))\,e\in \langle e\rangle\,and \,
[v-((v,e)/(e,e))\,e]\in \langle e\rangle^{\perp }.
\]
Thus $V= \langle e\rangle\oplus\langle e\rangle^{\perp }$ because
$\langle e\rangle\cap \langle e\rangle^{\perp }=(0)$. By the inductive
hypothesis applied to $\langle e\rangle^{\perp }$ which is of dimension $n-1$
there exist
 $e_{2}, e_{3},\dots , e_{n}\in \langle e\rangle^{\perp }$ which form
 a basis of $\langle e\rangle^{\perp }$ and such that
$(e_{i},e_{j})=0$ for $i,j=2,3, \dots ,n$\, and\,  $i\ne j$. Put $e_{1}=e$.
Then $e_{1},e_{2}, \dots ,e_{n}$ satisfy the required conditions.    
\end{Prop}
Put $k_{0}=\{\lambda\in k\mid \bar \lambda=\lambda \}$. Then $k_0$ is a
subfield of $k$ such that $dim k_{k_{0}}\leq 2$. For any $\lambda \in k$
define $\varphi(\lambda )=\lambda \bar \lambda$. Then
 $\varphi \colon\, k^*\to k_{0}^*$
and is a
homomorphism where for any field $K$, $K^*$ denotes the multiplicative
group of its non--zero elements. In case the involution of $k$ is not an
identity \ $\varphi(\lambda )=\prod_{\theta\in G(k/k_{0})}\theta(\lambda)$
where $G(k/k_{0})$ is the Galois group of the extension $k/k_{0}$.
So  $\varphi $ coincides with the ''usual'' norm homomorphism $N_{k/k_{0}}
 \colon\, k^*\to k_{0}^*$. 
\begin{Def}
Let $(\ ,\ ): V\times V\to k$ be a hermitian form. Then $v\in V$ is called
isotropic if \ $(v,v)=0$
\end{Def}
\begin{Prop}\label{P:b}
Let $(\ ,\ ): V\times V\to k$ be a hermitian form. If\\
$\varphi\colon\, k^*\to k_{0}^*$ is an epimorphism and dim $V\geq 2$ then $V$
contains a non--zero isotropic element.\\
Proof: Let $e_{1},e_{2},\dots , e_{n}\in V$ form a basis of \ $V$ and are
such that $(e_{i},e_{j})=0$ if $i\ne j$ (Proposition \ref{P:a}). Put
$\lambda_{i}=(e_{i},e_{i})$ for $i=1, 2,\ldots ,n$. We have $(e_{i},e_{i})=~
\overline {(e_{i},e_{i})}$. So $\lambda_{i}\in k_0$ for $i=1, 2,\ldots ,n$.
If $\lambda_{i}=0$
for some $i$ then the corresponding \ $e_{i}$ is isotropic. So we can assume
that $\lambda_{i}\ne 0$ for $i=1,2, \ldots ,n$. Let us consider
$v=xe_{i}+e_{j}$ where $x\in k$ such that $\varphi(x)=-\lambda_{j}/\lambda_{i}$
and $i\ne j$. Then $v\in V$ and $v\ne 0$. Moreover
\[
(v,v)=(xe_{i}+e_{j},xe_{i}+e_{j})=x\bar{x}\lambda_{i}+\lambda_{j}=
\varphi(x)\lambda_{i}+\lambda_{j}=\left(-\lambda_{j}/\lambda_{i}\right)\lambda_{i}+\lambda_{j}=0
\]
\end{Prop}
\begin{Def}
A field $K$ is called $C_{1}$ (\cite{S:1}) if every equation \\
$f(x_{1},x_{2},\ldots ,x_{n})=~0$ where $f$ is a homogenous polynomial of
degree $d<n$ with coefficients in $K$, admits a non-trivial solution in
$K^{n}$.
\end{Def}
\begin{Thm}
Let $(\ ,\ ): V\times V\to k$ be a hermitian form. If $k_{0}$ is $C_{1}$ and
dim~$V\geq 3$ then $V$ admits a non--zero isotropic element.\\
Proof:  $N_{k/k_{0}} \colon\, k^*\to k_{0}^*$ is an epimorphism since $k_{0}$
is $C_{1}$ (\cite{S:1}). It follows that
$\varphi\colon\, k^*\to k_{0}^*$ is an epimorphism
in case the involution is not an identity. So we can apply Propostion \ref{P:b}.
Suppose now that the involution is an identity and let
$e_{1},e_{2},\dots , e_{n}\in V$ form a basis of \ $V$ and are
such that $(e_{i},e_{j})=0$ if $i\ne j$ (Proposition \ref{P:a}). The equation
\[\lambda_{1}x_{1}^{2}+\lambda_{2}x_{2}^{2}+\cdots +\lambda_{n}x_{n}^{2}=0
\]
has a non--trivial solution $(x_{1},x_{2},\ldots ,x_{n})\, in\,  k=k_{0}$ where
$\lambda_{i}=(e_{i},e_{i})$ and $n\geq 3$. Let
$v=x_{1}e_{1}+x_{2}e_{2}+\cdots +x_{n}e_{n}\in V$. Then $v\ne 0$ and \\
 $(v,v)=
\lambda_{1}x_{1}^{2}+\lambda_{2}x_{2}^{2}+\cdots +\lambda_{n}x_{n}^{2}=0$.
\end{Thm}
\begin{Cor}
Let $k$ be a field of transcendence degree one over
an algebraically closed field $F$ with an involution which is an identity
on $F$. Suppose\\$(\ ,\ ):V\times V\to k$ is a hermitian form. If dim
$V\geq 3$ then $V$ admits a non--trivial isotropic element.\\
Proof: If $k$ is of transcendence degree one over $F$ then the same holds for
$k_{0}$ since $k$ is a finite extension of $k_{0}$.
By the Tsen Theorem (\cite{S:1}) $k_{0}$ is $C_{1}$, so the Theorem can be
applied.
\end{Cor}
\begin{Rem}
Suppose $k$ is a field with involution which satisfies the assumptions of
the Corollary and let $(\ ,\ ):V\times V\to k$ be a hermitian form on a vector
space $V$ with a trivial action of a finite group $G$.
Then obviously $W\subset W^{\perp }$ where $W$ denotes (a $G$--inwariant)
subspace generated by a non--zero isotropic element of \ $V$. So in case
dim $V\geq 3$ and the involution is an identity on $F$ (a condition very
natural if $k$ is considered a field of rational functions on a
non--singular curve over $F$) the proof of Maschke`s Theorem in case of
a real or complex field can not be generalized.
\end{Rem}

\end{document}